\documentclass{amsart}

\usepackage{graphicx}

\usepackage{amssymb}
\usepackage{amsfonts}
\usepackage[all]{xy}
\usepackage{amsmath}

\newtheorem{theorem}{Theorem}[section]
\newtheorem{lemma}[theorem]{Lemma}
\newtheorem{proposition}[theorem]{Proposition}
\newtheorem{corollary}[theorem]{Corollary}

\theoremstyle{definition}
\newtheorem{definition}[theorem]{Definition}
\newtheorem{example}[theorem]{Example}

\theoremstyle{remark}
\newtheorem{remark}[theorem]{Remark}

\numberwithin{equation}{section}

\newcommand{\K}{\mathbb{K}}
\newcommand{\NN}{\mathbb{N}}
\newcommand{\Z}{\mathbb{Z}}
\newcommand{\fr}{\mathfrak{F}}

\newcommand{\al}{\alpha}
\newcommand{\bt}{\beta}

\newcommand{\ups}{\upsilon}
\newcommand{\up}{\Upsilon}
\newcommand{\gm}{\gamma}

\newcommand{\A}{\mathcal{A}}

\newcommand{\C}{\mathcal{C}}

\newcommand{\M}{\mathcal{M}}

\newcommand{\nn}{\mathcal{N}}

\newcommand{\pd}[2]{#1^{(#2)}}

\begin{document}

\title{Generalized symmetric functions and invariants of matrices}


\author{Francesco Vaccarino}

\address{Politecnico di Torino - DISPEA, Corso Duca degli Abruzzi 24, I-10129, Torino, Italy}
              \email{francesco.vaccarino@polito.it}
\thanks{The author is supported by research grant Politecnico di Torino n.119 - 2004}

\date{}

\maketitle

\begin{abstract}
It is well known that over an infinite field the ring of symmetric
functions in a finite number of variables is isomorphic to the one
of polynomial functions on a single matrix that are invariants by
the action of conjugation by general linear group. We generalize
this result showing that the abelianization of the algebra of the
symmetric tensors of fixed order over a free associative algebra is
isomorphic to the algebra of the polynomials invariants of several
matrices over an infinite field or the integers. While proving the
main result we find generators and relations of abelianized divided
powers and symmetric products over any commutative ring.
\end{abstract}

\section{Introduction}\label{intro}
Let $\K$ be an infinite field and let $\K[y_1,\dots,y_n]^{S_n}$ be
the the ring of \emph{symmetric polynomials} in $n$ variables. The
general linear group $GL(n,\K)$ acts by conjugation on the full ring
$Mat(n,\K)$ of $n\times n$ matrices over $\K$. Denote by
$\K[Mat(n,\K)]^{GL(n,\K)}$ the ring of the polynomial invariants for
this actions. It is well known that
\begin{equation}\K[Mat(n,\K)]^{GL(n,\K)}\cong
\K[y_1,\dots,y_n]^{S_n}.\label{cla}\end{equation}

Let $M$ be a $\K-$module. Consider the tensor product
$T^n_{\K}(M):=M^{\otimes n}$, the symmetric group acts on
$T^n_{\K}(M)$ as a group of $\K-$module automorphisms and we denote
by $TS^n_{\K}(M):=(M^{\otimes n})^{S_n}$ the sub-$\K-$module of the
invariants for this action. The elements of $TS^n_{\K}(M)$ are
called the symmetric tensors of order $n$. If $M$ is a $\K-$algebra
then $TS^n_{\K}(M)$ is a sub$-\K-$algebra of $M$.

Let $F:=\K\{x_1,\dots,x_m\}$ be a free associative non commutative
algebra on $m$ variables. For a $\K-$algebra $A$ we write
$A^{ab}:=A/[A,A]$ for the abelianization of $A$, where $[A,A]$
denotes the ideal generated by the commutators. Consider
\begin{enumerate}
\item[(i)] $TS_{\K}^n(F)$ and
\item[(ii)]$TS_{\K}^n(F)^{ab}$ the abelianization of $TS^n_{\K}(F)$
\end{enumerate}
Note that when $m=1$ then $F$ is commutative and
$TS_{\K}^n(F)=TS_{\K}^n(F)^{ab}$.

We will find generators for both $TS_{\K}^n(F)$, $TS_{\K}^n(F)^{ab}$
and for divided powers of degree $n$ over any algebra. The result
presented are refinement of the one due to Ziplies \cite{zipgen} and
are proven in a rather simple way.

We also find relations, the first syzygy, for $TS_{\K}^n(F)^{ab}$
and then for the abelianization of divided powers of any algebra and
for divided powers of any commutative algebra.

Using the presentation we found we then prove the following
generalization of the isomorphism (\ref{cla}).
\begin{theorem}
Let $\K$ be an infinite field or the ring of integers and let the
general linear group $GL(n,\K)$ acts by simultaneous conjugation on
$m$ copies of $Mat(n,\K)$. Denote by $\K[Mat(n,\K)^m]^{GL(n,\K)}$
the ring of the invariants for this actions. Then
\[\K[Mat(n,\K)^m]^{GL(n,\K)}\cong TS^n_{\K}(\K\{x_1,\dots,x_m\})^{ab}\]
\end{theorem}

\begin{remark}
Let $\K$ be an infinite field, in {\cite{v3}} we proved that there
is an isomorphism $TS^n_{\K}(\K[x_1,\dots,x_m])\cong
\K[Z_{n,red}^m]^{GL(n,\K)}$ where $Z_{n,red}^m$ is the variety of
$m-$tuples of commuting $n\times n$ matrices. When $char\K=0$ then
we showed that the above isomorphism extends to the affine schemes
i.e. $TS^n_{\K}(\K[x_1,\dots,x_m])\cong \K[Z_n^m]^{GL(n,\K)}$, where
$Z_n^m$ is the affine scheme of $m-$tuples of commuting $n\times n$
matrices. Some links with the Hilbert scheme of points will be given
at ICM 2006 {\cite{v4}}
\end{remark}

\section{Conventions and notation}
All rings or algebras over a commutative ring will be with identity.
We denote by $\NN$ the natural numbers. Let $\K$ be a commutative
ring. We will use the following categories
\begin{itemize}
\item $\nn$, $\nn_{\K}$ the categories of non commutative rings and
non commutative $\K$-algebras
\item $\C$, $\C_{\K}$ the categories of commutative rings and commutative $\K$-algebras
\item $Mod_{\K}$ for the category of $\K-$modules
\item $Set$ the category of sets
\item We set $\A(B,C):=Hom_{\A}(B,C)$ in a category $\A$ with
$B,C\in Ob(\A)$ objects in $\A$.
\item For $A$ a set and any
additive monoid $M$, we denote by $M^{(A)}$ the  set of functions
$f:A\rightarrow M$ with finite support. \item Let $\al\in M^{(A)}$,
we denote by $\mid \al \mid$ the (finite) sum $\sum_{a\in A}
\al(a)$.
\end{itemize}
\begin{remark}
In order to avoid cumbersome notation we limit ourself to finitely
generated rings and algebras, but we remark that our results remain
valid over not finitely generated rings and algebras, unless
otherwise stated.
\end{remark}

\section{Symmetric functions}
Set $\Lambda_n:=\K[y_1,\dots,y_n]^{S_n}$ for the ring of
\emph{symmetric polynomials} in $n$ variables: it is freely
generated by the elementary symmetric functions $e_1,\dots,e_n$
given by the equality
\begin{equation}\label{elem}
\sum_{k=0}^n t^k e_k:=\prod_{i=1}^n(1+ty_i)
\end{equation}
Here $e_0=1$ and $t$ is an independent variable (see {\cite{M}}).
Furthermore one has
\begin{equation}\label{elem2}
e_k(y_1,\dots,y_n)=\sum_{i_1<i_2<\dots <i_k\leq
n}y_{i_1}y_{i_2}\cdots y_{i_k}
\end{equation}

The action of $S_n$ on $\K[y_1,y_2,\dots ,y_n]$ preserves the usual
degree. We denote by $\Lambda_{n}^k$ the $\K$-submodule of
invariants of degree $k$.

Let $q_n:\K[y_1,y_2,\dots ,y_n]\rightarrow \K[y_1,y_2,\dots
,y_{n-1}]$ be given by
\[\begin{cases} y_n\mapsto 0
\\y_i\mapsto y_i {\hbox{ for }}i=1,\dots,n-1 \end{cases}\]

for $i=1,\dots,n-1$. One has $q_n(\Lambda_{n}^k)=\Lambda_{n-1}^k$
and it is easy to see that $\Lambda_{n}^k\cong\Lambda_{k}^k$ for all
$n \geq k$. Denote by $\Lambda^k$ the limit of the inverse system
obtained in this way.

\begin{definition}The ring $\Lambda:=\bigoplus_{k\geq 0}\Lambda^k$ is called the
ring of symmetric functions (over $\K$)
\end{definition}

\noindent It can be shown {\cite{M}} that $\Lambda_R$ is a
polynomial ring, freely generated by the (limits of the) $e_k$, that
are given by
\begin{equation} \sum_{k=0}^{\infty} t^k e_k:=\prod_{i=1}^{\infty}(1+tx_i).
\end{equation}
Furthermore the kernel of the projection $\pi_n:\Lambda \rightarrow
\Lambda_{n}$ is generated by the $e_{n+k}$, where $k\geq 1$.

We have another distinguished kind of functions in $\Lambda_n$
beside the elementary symmetric ones: the \emph{power sums}, for any
$r\in \NN$ the $r$-th power sum is
\begin{equation}\label{pow}
p_r:=\sum_{i\geq 1}y_i^r\end{equation} Let $g\in \Lambda_n$, set $g
\cdot p_r=g(y_1^r,y_2^r,\dots,y_n^r)$, this is again a symmetric
function. Since the $e_i$ generate $\Lambda_n$ we have that $g\cdot
p_r$ can be expressed as a polynomial in the $e_i$ we set
\begin{equation} P_{h,k}:=e_h\cdot p_k\label{plet}\end{equation}
to denote it.

The monomials form a $\K-$basis of $\K[y_1,\dots,y_n]$ permuted by
$S_n$. Thus, the sums of monomials over the orbits form a $\K-$basis
of the ring $\Lambda_n$ and their limits form a basis of $\Lambda$.
Let $y_1^{\lambda_1}y_2^{\lambda_2}\cdots y_n^{\lambda_n}$ be a
monomial, after a suitable permutation we can suppose $\lambda_1\geq
\lambda_2\geq\dots\geq\lambda_n\geq 0$. We set $m_{\lambda}$ for the
orbit sum corresponding to such
$\lambda=(\lambda_1,\lambda_2,\dots,\lambda_n)\in\NN^n$ then
\begin{equation}\mathcal{P}_n=\{m_{\lambda}\,:\,\lambda_1\geq \lambda_2\geq\dots\geq\lambda_n\geq 0,\,\lambda_i\in\NN\}
\end{equation} and
\begin{equation}\mathcal{P}_{n,k}=\{m_{\lambda}\,:\,\sum_i\lambda_i=k\}\end{equation}
are $\K-$bases of $\Lambda_n$ and $\Lambda_{n}^k$ respectively. As
before the limits of the $m_{\lambda}$ form a basis of $\Lambda$ and
$\Lambda^k$ and $\ker(\Lambda\rightarrow\Lambda_n)$ has basis
$\{m_{\lambda}\,:\,\exists k>n{\hbox{ with }} \lambda_k>0\}$

\section{Symmetric Tensors on Free Algebras}
We set here a generalization of $\Lambda$ and $\Lambda_n$, our
exposition will be streamlined on the one given in the previous
section.

Let $F:=\K\{x_1,\dots,x_m\}$ be a free associative non commutative
algebra. Let $k\in \NN$, we denote by $\mathbf{f}$ the sequence
$(f_1\dots,f_k)$ in $F$ and by $\al$ the element
$(\al_1,\dots,\al_k)\in\NN^k$, with $\mid\al\mid:=\sum\al_j \leq n$
. Let $t_1,\dots,t_k$ be commuting independent variables, we set as
usual $t^{\al}:=\prod_{i}t_i^{\al_i}$. We define elements
$e_{\al}^n(\mathbf{f})\in TS_{\K}^n(F)$ by
\begin{equation}\sum_{\al}t^{\al}\otimes e_{\al}^n(\bold{f}) :=
(1+\sum_ht_h\otimes f_h)^{\otimes n}\end{equation}
\begin{example}
Let $f,g\in F$ then
\begin{eqnarray*}
e_{(2,1)}^3(f,g)&=& f\otimes f\otimes g+f\otimes g\otimes
f+g\otimes f\otimes f \\
e_{(2,1)}^4(f,g)&=& f\otimes f\otimes g\otimes 1+f\otimes g\otimes
f\otimes 1+g\otimes f\otimes f\otimes
1\\
& + &f\otimes f\otimes 1\otimes g+f\otimes g\otimes 1\otimes
f+g\otimes f\otimes 1\otimes f\\
&+&f\otimes 1\otimes f\otimes g+f\otimes 1\otimes g\otimes
f+g\otimes 1\otimes f\otimes f\\
&+&1\otimes f\otimes f\otimes g + 1\otimes f\otimes g\otimes
f+1\otimes g\otimes f\otimes f\\
\end{eqnarray*}
\end{example}
\begin{lemma} The element
$e^n_{(\al_1,\dots,\al_k)}(f_1,\dots,f_k)$ is the orbit sum (under
the considered action of $S_n$) of
\[f_1^{\otimes\al_1}\otimes f_2^{\otimes\al_2}\cdots \otimes f_k^{\otimes\al_k}\otimes
1^{\otimes ( n-\sum_i\al_i)}\]
\end{lemma}\label{orbit}
\begin{proof} Let $E$ be the set of mappings
$\phi:\left\{1,\dots,n\right\}\rightarrow\left\{1,\dots,k+1\right\}$.
We define a mapping $\phi\mapsto\phi^*$ of $E$ into $\NN^{k+1}$ by
putting $\phi^*(i)$ equal to the cardinality of $\phi^{-1}(i)$. For
two elements $\phi_1,\phi_2$ of $E$, to satisfy $\phi_1^*=\phi_2^*$
it is necessary and sufficient that there should exist $\sigma\in
S_n$ such that $\phi_2=\phi_1\circ\sigma$. Set $f_{k+1}:=1_R$ and
$E(\al):=\left\{\phi\in E \mid \phi^*=(\al_1,\dots,\al_k,
n-\sum_i\al_i)\right\}$, then we have
\[e_{\alpha}(\bold{f}) =\sum_{\phi\in E(\al)}f_{\phi(1)}\otimes
f_{\phi(2)}\otimes\cdots \otimes f_{\phi(n)}\] and the lemma is
proved.\qed
\end{proof}
\begin{definition}\label{degree}
Let $\Upsilon$ denote the set of monomials in $F$. For $\upsilon\in
\Upsilon$ let $\partial_i(\ups)$ denote the degree of $\ups$ in
$x_i$, for all $i=1,\dots ,m$. We set
\begin{equation}
\partial(\ups):=(\partial_1(\ups),\dots ,\partial_m(\ups))
\end{equation}
for its \emph{multidegree}. The total degree of $\ups$ is
$\sum_i\partial_i(\ups)$. We denote by $\Upsilon^+$ be the set of
monomials of positive degree.
\end{definition}
Is is clear that $\up$ is a $\K-$basis of $F$ so that
$\up_n:=\{\ups_1\otimes \ups_2\otimes\cdots \otimes
\ups_n\,:\,\ups_j\in\up \}$ is a $\K-$basis of $T^n_{\K}(F)$
permuted by $S_n$. Thus, the sums of the elements of $\up_n$ over
their orbits form a $\K-$basis of the symmetric tensors
$TS^n_{\K}(F)$.

Let $\al\in \NN^{(\up^+)}$, then there exist $k\in \NN$ and
$\ups_1,\dots,\ups_k\in\up^+$ such that $\al(\ups_i)=\al_i\neq 0$
for $i=1,\dots,k$ and $\al(\up)=0$ when $\ups\neq
\ups_1,\dots,\ups_k$. We set
\begin{equation}\label{coef}e_{\al}^n:=e_{(\al_1,\dots,\al_k)}^n(\ups_1,\dots,\ups_k),\end{equation}
then
\begin{equation}\sum_{\mid\al\mid\leq n}t^{\al}\otimes e_{\al}^n
=(1+\sum_{\ups\in\up^+}t_{\ups}\otimes\ups)^{\otimes
n}\end{equation} where $t_{\ups}$ are commuting independent
variables indexed by monomials and
\begin{equation}t^{\al}:=\prod_{\ups\in\up^+}t_{\ups}^{\al(\ups)}\end{equation}
for all $\al\in \NN^{(\up^+)}$.
\begin{proposition} The set
\[\mathcal{B}_n:=\{e_{\al}^n\;:\;\mid \al\mid\leq n\}\] is a
$\K-$basis of $TS_{\K}^n(F)$.
\end{proposition}
\begin{proof}
By Lemma\,\ref{orbit} and (\ref{coef}), the $e_{\al}^n$ are a
complete system of representatives (for the action of $S_n$) of the
orbit sums of the elements of $\up_n$.\qed
\end{proof}
\begin{definition}
We define a $S_n-$invariant multidegree on $T^n(F)$ by
\[\partial(1^{\otimes i}\otimes x_j\otimes
1^{n-i-1}):=\partial(x_j)\] for all $i,j$. We denote by
$T^n_{\K}(F)_{\delta}$ (resp. $TS^n_{\K}(F)_{\delta}$) the linear
span of the elements of multidegree $\delta\in\NN^m$.
\end{definition}
\begin{proposition}
The set
\[\mathcal{B}_{n,\delta}:=\{e_{\al}^n\;:\;\mid \al\mid\leq n\text{ and } \partial(e_{\al})=\delta\}\] is a
$\K-$basis of $TS_{\K}^n(F)_{\delta}$ for all $\delta\in\NN^m$.
\end{proposition}
\begin{proof}
Observe that
$\partial(e_{\al})=\sum_{\ups\in\up^+}\al_{\ups}\partial(\ups)$ and
apply the above Proposition.\qed
\end{proof}
\section{Generators}
First of all we compute the product between the elements of
$\mathcal{B}_n$
\begin{proposition}[Product Formula]
Let $h,k\in\NN$, $\al\in\NN^h$, $\bt\in\NN^k$ be such that
$\mid\al\mid,\,\mid\bt\mid\leq n$. Let
$r_1,\dots,r_h,s_1,\dots,s_k\in F$. Set again
\[e_{\al}^n(\mathbf{r}):=e^n_{(\al_1,\dots,\al_h)}(r_1,\dots,r_h)
\text{ and }
e_{\bt}^n(\mathbf{s}):=e^n_{(\bt_1,\dots,\bt_k)}(s_1,\dots,s_k)\]
then \[e_{\al}^n(\mathbf{r})e_{\bt}^n(\mathbf{s})=\sum_{\gm}
e_{\gm}(\mathbf{r},\mathbf{s},\mathbf{rs})\] where
\begin{eqnarray*}\mathbf{rs} &:=&
(r_1s_1,r_1s_2,\dots,r_1s_k,r_2s_1,\dots,r_2s_k,\dots,r_hs_k)\\
\gamma &:=&
(\gamma_{10},\dots,\gamma_{h0},\gamma_{01},\dots,\gamma_{0k},
\gamma_{11},\gamma_{12},\dots,\gamma_{hk})
\end{eqnarray*}
 are such that
\begin{equation}\label{sys}
\begin{cases} \gamma_{ij}\in \NN \\
\mid \gamma \mid \leq n \\
\sum_{j=0}^k \gamma_{ij}=\al_i \;\;{\text{for}}\;\; i=1,\dots,h\\
\sum_{i=0}^h \gamma_{ij}=\bt_j \;\; {\text{for}}\;\; j=1,\dots,k.
\end{cases}
\end{equation}
\end{proposition}
\begin{proof}
Let $t_1,t_2$ be two commuting independent variables and let $a,b\in
F$ then one has
\begin{equation}\label{aperb}(1+t_1\otimes a)^{\otimes
n}(1+t_2\otimes b)^{\otimes n}=(1+t_1\otimes a+t_2\otimes
b+t_1t_2\otimes ab)^{\otimes n}\end{equation} hence
\begin{align*}
(1+\sum_{i=1}^n t_1^i\otimes e_i^n(a))(1+\sum_{j=1}^n t_2^j &\otimes e_j^n(b))\\
&=1+\sum_{i,j}t_1^it_2^j\otimes e_i^n(a)e_j^n(b)\\
&=1+\sum _{l_1,l_2,l_{12}}t_1^{l_1+l_{12}}t_2^{l_2+l_{12}}\otimes
e^n_{_{(l_1,l_2,l_{12})}}(a,b,ab)
\end{align*}
The desired equation the easily follows.\qed
\end{proof}
\begin{remark}
The above Product Formula can be derived from the one of N.Roby
{\cite{rs}} in the context of divided powers. It has been also
derived by Ziplies in his papers \cite{zipgen} on the divided powers
algebra $\widehat{\Gamma}$.
\end{remark}
\begin{corollary} Let $k\in \NN$, $a_1,\dots ,a_k\in F$,
$\al=(\al_1,\dots,\al_k)\in \NN^k$ with $\mid\al\mid\leq n$. Then
$e^n_{(\al_1,\dots,\al_k)}(a_1,\dots,a_k)$ belongs to the subalgebra
of $TS^n_{\K}(F)$ generated by the $e_i^n(\ups)$, where
$i=1,\dots,n$ and $\ups$ is a monomial in the $a_1,\dots,a_k$.
\end{corollary}
\begin{proof}
We prove the claim by induction on $\mid\al\mid$ assuming that
$\al_i>0$ for all $i$. (note that $1 \leq k\leq \sum_j\al_j$). Since
$n$ is fixed we suppress the superscript $n$ for all the proof.

If $\sum_j\al_j=1$ then $k=1$ and
$e_{(\al_1,\dots,\al_k)}(a_1,\dots,a_k)=e_1(a_1)$. Suppose the claim
true for all $e_{(\bt_1,\dots,\bt_h)}(b_1,\dots,b_h)$ with
$b_1,\dots,b_h\in F$ and $\mid\bt\mid < \mid\al\mid$. Let $k,
a_1,\dots ,a_k, \al$ be as in the statement, then we have by the
Product Formula
$$e_{\al_1}(a_1)e_{(\al_2,\dots,\al_k)}(a_2,\dots,a_k)=$$
$$=e_{(\al_1,\dots,\al_k)}(a_1,\dots,a_k)+\sum
e_{\gamma}(a_1,\dots,a_k,a_1a_2,\dots,a_1a_k),$$ where
$$\gamma=(\gamma_{10},\gamma_{01},\dots,\gamma_{0h},\gamma_{11},\gamma_{12},\dots,\gamma_{1h})$$
with $h=k-1$, $\sum_{j=0}^{h} \gamma_{1j}=\al_1$ with $\sum_{j=1}^h
\gamma_{1j}>0$, and $\gamma_{0j}+\gamma_{1j}=\al_j$ for
$j=1,\dots,h$. Thus
$$\gamma_{10}+\gamma_{01}+\dots+\gamma_{0h}+\gamma_{11}+\dots+\gamma_{1h}
=\sum_j\al_j - \sum_{j=1}^h \gamma_{1j}<\sum_j\al_j.$$ Hence
$$e_{(\al_1,\dots,\al_k)}(a_1,\dots,a_k)=$$
$$e_{\al_1}(a_1)e_{(\al_2,\dots,\al_k)}(a_2,\dots,a_k)
-\sum e_{\gamma}(a_1,\dots,a_k,a_1a_2,a_1a_3,\dots,a_1a_k),$$ where
$\mid\gm\mid=\sum_{r,s}\gm_{rs}<\mid\al\mid$. So the claim follows
by induction hypothesis.\qed
\end{proof}
\begin{corollary}\label{gen}
The algebra of symmetric tensors $TS^n(F)$ of order $n$ is generated
by the $e_i^n(\ups)$ where $1\leq i\leq n$ and $\ups\in\up^+$
\end{corollary}
\begin{proof}
It follows from the above corollary applied to the elements of the
basis $\mathcal{B}_n$.\qed
\end{proof}
\begin{remark}
The above corollaries can be proved using Cor(4.1) and Cor(4.5) in
\cite{zipgen}
\end{remark}
\begin{lemma}\label{ple} For all $f\in F$, and $k,h\in \NN$,
$e_h^n(f^k)$ belongs to the subalgebra of $TS^n_{\K}(F)$ generated
by the $e_j^n(f)$.
\end{lemma}
\begin{proof}
Let $f\in F$ there is an obvious $S_n-$equivariant morphism
\begin{equation}\label{eval}\rho_f:\begin{cases}
    \K[y_1,\dots,y_m]\rightarrow T^n_{\K}F,   \\
    y_h\mapsto 1^{\otimes(h-1)}\otimes f \otimes 1^{\otimes(n-h)}
\end{cases}
\end{equation}
induced by the evaluation $\K[x]\rightarrow \K[f]$. It is the clear
that $\rho_f(e_i)=e_i^n(f)$ so that
\[e_h^n(f^k)=\rho_f(e_h\cdot p_k)=\rho_f(P_{h,k})\]
and the result is proved.\qed\end{proof}
\begin{definition}
A monomial $\ups\in\up^+$ is called \emph{primitive} if it is not
the power of another one.
\end{definition}
\begin{example}$x_1x_2x_1x_2$ is not primitive while $x_1x_2x_1x_1$ is
primitive.
\end{example}
We have then the following refinement of Cor.\ref{gen}.
\begin{theorem}[Generators]\label{prigen}
The algebra $TS^n_{\K}(F)$ is generated by $e_i^n(\ups)$ with $1\leq
i \leq n$ and $\ups$ primitive
\end{theorem}
\begin{proof}
It follows from Cor.\ref{gen} and Lemma\ref{ple}.\qed
\end{proof}
\subsection{Abelianization} Given any  $\K-$algebra $R$ one can form
the bilateral ideal $[R,R]$ generated by the commutators
$[a,b]:=ab-ba$, $a,b\in R$.

We set $R^{ab}:=R/[R,R]$ and call it \emph{the abelianization} of
$R$. The abelianization of $R$ is commutative and together with the
homomorphism $\mathfrak{ab}:R\longrightarrow R/[R,R]$ has the
universal property that the map $\C_{\K}(R^{ab},S)\rightarrow
\nn_{\K}(R,S)$ given by $\varphi\mapsto\varphi\circ\mathfrak{ab}$ is
an isomorphism for all $S\in\C_{\K}$.

\begin{definition}
Consider $\up/\sim$ the set of the equivalence classes of monomials
$\ups\in\up^+$ where $\ups\sim\ups'$ if and only if there a cyclic
permutation $\sigma$ such that $\sigma(\ups)=\ups'$. We set $\Psi$
to denote the set of equivalence classes in $\up^+/\sim$ made of
\emph{primitive} monomials
\end{definition}
We are now able to give the main result of this paragraph.
\begin{theorem}\label{genab}
The algebra $(TS^n_{\K}(F))^{ab}$ is generated by $e_i^n(\ups)$ with
$1\leq i \leq n$ and $\ups$ that varies in a complete set of
representatives of $\Psi$.
\end{theorem}
\begin{proof}
Let $\mathfrak{ab}:TS_{\K}(F)\longrightarrow (TS^n_{\K}(F))^{ab}$ be
the canonical projection, using (\ref{aperb}) it easy to see that
\begin{equation}\mathfrak{ab}(e_i^n(rs))=\mathfrak{ab}(e_i^n(sr))\end{equation}
for all $1\leq i \leq n$ and $r,s\in F$. The result then follows
from Th.\ref{prigen} and Prop.\ref{loth} by the surjectivity of
$\mathfrak{ab}$.\qed
\end{proof}

\section{Polynomial Laws}\label{poly}
To link symmetric tensors to linear representations we shall use the
determinant so that we are lead to the topic of polynomial laws: we
recall the definition of this kind of map between $\K$-modules that
generalizes the usual polynomial mapping between free $\K$-modules
(see {\cite{bo,rl,rs}}).
\begin{definition} Let $A$ and $B$ be two
$\K$-modules. A \emph{polynomial law} $\varphi$ from $A$ to $B$ is a
family of mappings $\varphi_{_{L}}:L\otimes_{\K} A \longrightarrow
L\otimes_{\K} B$, with $L\in\C_{\K}$ such that the following diagram
commutes
\begin{equation}
\xymatrix{
  L\otimes_{\K}A \ar[d]_{f\otimes id_A} \ar[r]^{\varphi_L}
                & L\otimes_{\K} B \ar[d]^{f\otimes id_B}  \\
  M\otimes_{\K} A \ar[r]_{\varphi_M}
                & M\otimes_{\K} B             }
\end{equation}
for all $L,\,M\in\C(\K)$ and all $f\in \C_{\K}(L,M)$.
\end{definition}
\begin{definition}
Let $n\in \NN$, if $\varphi_L(au)=a^n\varphi_L(u)$, for all $a\in
L$, $u\in L\otimes_{\K} A$ and all $L\in\C_{\K}$, then $\varphi$
will be said \emph{homogeneous of degree} $n$.
\end{definition}
\begin{definition}
If $A$ and $B$ are two $\K$-algebras  and
\[
\begin{cases} \varphi_L(xy)&=\varphi_L(x)\varphi_L(y)\\
         \varphi_L(1_{L\otimes A})&=1_{L\otimes B}
         \end{cases}
         \]
for $L\in\C_{\K}$ and for all $x,y\in L\otimes_{\K} A$, then
$\varphi$ is called \emph{multiplicative}.
\end{definition}
\begin{remark} A polynomial law $\varphi:A\rightarrow B$  is a
natural transformation $-\otimes_{\K}A\rightarrow -\otimes_{\K}B$.
\end{remark}
Let $A$ and $B$ be two $\K$-modules and $\varphi:A\rightarrow B$ be
a polynomial law. The following result on polynomial laws is a
restatement of Th\'eor\`eme I.1 of {\cite{rl}}.
\begin{theorem}\label{roby} Let $S$ be a set.
\begin{enumerate}
\item Let $L=\K[x_s]_{s\in S}$ and let $\{a_{s}\, :\,s\in S\}\subset A$
be such that $a_{s}=0$ except for a finite number of $s\in S$, then
there exist $\varphi_{\xi}((a_{s}))\in B$, with $\xi \in \NN^{(S)}$,
such that:
\[\varphi_{_{L}}(\sum_{s\in S} x_s\otimes
a_{s})=\sum_{\xi \in \NN^{(S)}}  x^{\xi}\otimes
\varphi_{\xi}((a_{s})),\] where $x^{\xi}:=\prod_{s\in S}
x_s^{\xi_s}$.
\item Let $R$ be any commutative $\K$-algebra and let
$(r_s)_{s\in S}\subset R$, then: \[\varphi_{_{R}}(\sum_{s\in S}
r_s\otimes a_{s})=\sum_{\xi \in \NN^{(S)}}  r^{\xi}\otimes
\varphi_{\xi}((a_{s})),\] where $r^{\xi}:=\prod_{s\in S}
r_s^{\xi_s}$.
\item If $\varphi$ is homogeneous of degree $n$, then one has $\varphi_{\xi}((a_{s}))=0$ if $\mid \xi \mid$ is
different from $n$. That is: \[\varphi_{_{R}}(\sum_{a\in A}
r_a\otimes a)=\sum_{\xi \in \NN^{(A)},\,\mid \xi \mid=n}
r^{\xi}\otimes \varphi_{\xi}((a)).\] In particular, if $\varphi$ is
homogeneous of degree $0$ or $1$, then it is constant or linear,
respectively.
\end{enumerate}
\end{theorem}
\begin{remark}
The above theorem means that a polynomial law $\varphi:A\rightarrow
B$ is completely determined by its coefficients
$\varphi_{\xi}((a_{s}))$, with $(a_s)_{s\in S} \in S^{(\NN)}$.
\end{remark}
\begin{remark} If $A$ is a free $\K$-module and $\{a_{t}\,
:\, t\in T\}$ is a basis of $A$, then $\varphi$ is completely
determined by its coefficients $\varphi_{\xi}((a_{t}))$, with $\xi
\in \NN^{(T)}$. If also $B$ is a free $\K$-module with basis
$\{b_{u}\, :\, u\in U\}$, then $\varphi_{\xi}((a_{t}))=\sum_{u\in
U}\lambda_{u}(\xi)b_{u}$. Let $a=\sum_{t\in T}\mu_{t}a_{t}\in A$.
Since only a finite number of $\mu_{t}$ and $\lambda_{u}(\xi)$ are
different from zero, the following makes sense:
\begin{eqnarray*}\varphi(a)=\varphi(\sum_{t\in T}\mu_{t}a_{t}) =  \sum_{\xi\in
\NN^{(T)}} \mu^{\xi}\varphi_{\xi}((a_{t}))& = & \sum_{\xi\in
\NN^{(T)}} \mu^{\xi}(\sum_{u\in U}\lambda_{u}(\xi)b_{u})\\ & = &
\sum_{u\in U}(\sum_{\xi\in \NN^{(T)}}\lambda_{u}(\xi)
\mu^{\xi})b_{u}.\end{eqnarray*} Hence, if both $A$ and $B$ are free
$\K$-modules, a polynomial law $\varphi:A\rightarrow B$ is simply a
polynomial map.
\end{remark}
\begin{definition}\label{polset}
Let $\K$ be a commutative ring.
\begin{itemize}
\item[(1)] For $M,N$ two $\K-$modules we set $P^n(M,N)$ for
the set of homogeneous polynomial laws $M\rightarrow N$ of degree
$n$
\item[(2)] If $M,N$ are two $\K-$algebras we set $MP^n(M,N)$ for
the multiplicative homogeneous polynomial laws $M\rightarrow N$ of
degree $n$.
\end{itemize}
\end{definition}
The assignment $N\longrightarrow P^n(M,N)$ (resp. $N\longrightarrow
MP^n(M,N)$) determines a functor from $\K-$modules (resp.
$\K-$algebras) to $Set$.
\section{Divided Powers}
The functors just introduced in Def.\ref{polset} are represented by
the divided powers which we introduce right now in the stream.
\begin{definition}
Let $\K$ be any commutative ring with identity. For a $\K$-module
$M$ let $\Gamma_{\K}(M)$ denote its divided powers algebra (see
{\cite{bo}}, {\cite{rl}} and {\cite{rs}}). $\Gamma_{\K}(M)$ is an
associative and commutative $\K$-algebra with identity $1_{\K}$ and
product $\times$, with generators $\pd m k $, with $m\in M$, $k \in
\Z$ and relations, for all $m,n\in M$:
\begin{enumerate}
\item $\pd m i  = 0, \forall i<0$;
\item $\pd m 0  = 1_{\scriptstyle{\K}},  \forall m\in M$;
\item $\pd {(\al m)} i  = \al^i\pd m i,  \forall \al\in \K,
\forall i\in \NN$;
\item $\pd {(m+n)} k  = \sum_{i+j=k}\pd m i \pd n j ,  \forall
k\in \NN$;
\item $\pd m i \times \pd m j  = {i+j\choose i}\pd m {i+j} ,
 \forall i,j\in \NN$.
\end{enumerate}
\end{definition}
\noindent The $\K$-module $\Gamma_{\K}(M)$ is generated by finite
products $\times_{i\in I} \pd {x_i} {\al_i}$ of the above
generators. The divided powers algebra $\Gamma_{\K}(M)$ is a
$\NN$-graded algebra with homogeneous components $\Gamma^k_{\K}(M)$,
($k\in \NN$), the submodule generated by $\{\times_{i\in I} \pd
{x_i} {\al_i}\; :\;\mid\al\mid=\sum_i\al_i=k\}$. One easily check
that $\Gamma_{\K}$ is a functor from $\K$-modules to commutative
graded $\K$-algebras.
\subsection{Universal properties}
The following properties give the motivation for the introduction of
divided powers in our setting.
\subsubsection{}
$\Gamma_{\K}^n$ is a covariant functor from $Mod_{\K}$ to $Mod_{\K}$
and one can easily check that it preserves surjections. The map
$\gamma^n:r\mapsto \pd r n$ is a polynomial law
$M\rightarrow\Gamma^n_{\K}(M)$ homogeneous of degree $n$: we call it
\emph{the universal one}. We have an isomorphisms
$Mod_{\K}(\Gamma^n_{\K}(M),N)\rightarrow P^n_{\K}(M,N)$ given by
$\phi\mapsto\phi\circ\gamma^n$.
\subsubsection{}
If $R$ is $\K-$algebra then $\Gamma^n_{\K}(R)$ inherits a structure
of $\K-$algebra by induced by setting $\pd a n\pd b n = \pd{(ab)}n$.
The unit in $\Gamma_{\K}^n(R)$ is $\pd 1 n$. It was proven by N.Roby
{\cite{rs}} that in this way $R\to\Gamma_{\K}^n(R)$ gives a functor
from $\K$-algebras to $\K-$algebras such that
$\gamma^n(a)\gamma^n(b)=\gamma^n(ab)$, $\forall a,b \in R$. Hence
$\gamma^n\in MP^n_{\K}(M,\Gamma^n_{\K}(M))$ and the map
$\nn_{\K}(\Gamma_{\K}^n(R),S)\rightarrow MP_{\K}(R,S)$ given by
$\phi\mapsto\phi\circ\gamma^n$ is an isomorphism for all $R,S\in
\nn_{\K}$.
\subsection{Flatness and Symmetric Tensors}
Let $M\in Mod_{\K}$. Consider again the symmetric tensors
$TS^n_{\K}(M)$ of order $n$.  When $M\in\nn_{\K}$ (resp.
$M\in\C_{\K}$) then $TS^n_{\K}(M)\in\nn_{\K}$ (resp.
$TS^n_{\K}(M)\in\C_{\K}$). The homogeneous polynomial law
$M\rightarrow TS^n_{\K}(M)$ given by $x\mapsto x^{\otimes n}$ gives
a morphism $\tau_n:\Gamma_{\K}^n(M)\rightarrow TS^n_{\K}(M)$ that is
an isomorphism when $M$ is flat over $\K$.(see {\cite{de}} 5.5.2.5
pag.123). Note that $\tau_n(\pd r n)=r^{\otimes n}$ (see
{\cite{bo}}) and observe that this implies that $\tau_n(1+t\otimes
r)^{(n)}=(1+t\otimes r)^{\otimes n}$ where $t$ is a commuting
independent variable so that
\begin{equation}\label{egia}
\tau_n(\pd 1 {n-j}\pd {r} j)=e_i^n(r)
\end{equation}
for all $\K-$algebra $R$, $\forall r\in R$ and $j=1,\dots,n$.
\begin{corollary}\label{sip}
Let $R$ be a $\K-$algebra generated by $\{r_i\}_{i\in I}$ then
\begin{enumerate}
\item $\Gamma_{\K}^n(R)$ is generated by $\pd 1 {n-j}\pd {\ups} j$ where
$1\leq j\leq n$ and $\ups$ varies in the set of primitive monomials
in $\{r_i\}_{i\in I}$
\item $(\Gamma_{\K}^n(R))^{ab}$ is
generated by $\pd 1 {n-j}\pd {\ups} j$ where $\ups$ varies in a
complete set of representative of equivalence classes (under cyclic
permutations) of primitive monomials in the $\{r_i\}_{i\in I}$.
\end{enumerate}
\end{corollary}
\begin{proof}
Given a $\K-$algebra $R$ and a set $\{r_i\}_{i\in I}$ of its
generators we have a surjective homomorphism
$F_I:=\K\{x_{\al}\}\twoheadrightarrow R$ hence another one
$TS^n_{\K}(F_I)\cong\Gamma_{\K}^n(F_I)\twoheadrightarrow\Gamma_{\K}^n(R)$
and the result follows from Cor.\ref{prigen} recalling (\ref{egia}).
The result on the abelianization follows from Th.\ref{genab} since
the abelianization functor preserves surjections.\qed
\end{proof}
\begin{remark}
Corollary \ref{sip}(1) is a refinement of Cor.(4.5) in
\cite{zipgen}.
\end{remark}
\begin{corollary}\label{rh}Let $R$ be commutative.
\begin{enumerate}
\item
The ring
$TS_{\K}^n(\K[y_1,\dots,y_m])\cong\K[x_{11},x_{12},\dots,x_{m1},\dots,x_{mn}]^{S_n}$
of the multisymmetric functions also known as the ring of the vector
invariants of the symmetric group is generated by the $e_i^n(\ups)$
with $\ups=y_1^{\al_1}\cdots y_{m}^{\al_m}$ such that
$\al_1,\dots,\al_m$ are coprime.
\item Let $R$ be a commutative $\K-$algebra generated by
$\{r_i\}_{i\in I}$ then $\Gamma_{\K}^n(R)$ is generated by $\pd 1
{n-j}\pd {\ups} j$ where $\ups=\prod_ir_i^{\al_i}$ is such that
$\sum_i\al_i$ is finite and the $\al_i$ are coprime.
\end{enumerate}
\end{corollary}
\begin{proof}
It follows from Th.\ref{genab} since $\Gamma^n$ preserves
surjections. Write $P=K[y_1,\dots,y_m]$, then $P$ is a free
$\K-$module hence the surjection $F\twoheadrightarrow P$ induces
another surjection $TS^nF^{ab}\twoheadrightarrow TS^nP$. Use
Th.\ref{genab}.\qed
\end{proof}
\begin{remark}
This gives another proof of Th.1 in \cite{vf}
\end{remark}
\begin{remark}
Corollary \ref{rh} has been improved by D.Rydh \cite{rd} which found
a {\em minimal\/} generating set for $TS_{\K}^n(\K[y_1,\dots,y_m])$.
\end{remark}
\subsection{Good Characteristics}
In the ring $\Lambda$ of the symmetric functions it holds the
following well known \textbf{Newton's Formulas}
\begin{equation}\label{newt}
(-1)^kp_{k+1} + \sum_{i=1}^k(-1)^ip_ie_{k+1-i}
=(k+1)e_{k+1}\end{equation} for all $k>0$. It is clear that these
equalities hold also in $\Lambda_n$ with $e_i=0$ for $i>n$.
\begin{proposition}
Let $R$ be generated by $\{r_i\}_{i\in I}$ as a $\K-$algebra. Let
$n!$ be invertible in $\K$ then $\Gamma_{\K}^n(R)$ is generated by
$\pd 1 {n-1}\pd{\ups}1$ where $\ups$ varies in the monomials in the
$r_i$. If $R$ is commutative the monomial $\ups$ can be taken of
degree not greater than $n$.
\end{proposition}
\begin{proof}
In this case we have that Newton's formulas made $p_1,p_2,\dots,p_n$
a generating set for $\Lambda_n$ hence $\pd 1{n-i}\pd {\ups}i$
belongs to the subring generated by the $\pd 1 {n-1}\pd {(\ups^k)}1$
and this gives the desired result when combined with the results on
generating sets given in the previous paragraphs. If $R$ is
commutative the result follows because it holds on
$\Gamma_{\K}^n(\K[y_1,\dots,y_m])$ by E.Noether's bound.\qed
\end{proof}
\section{$\Gamma^n$ vs $TS^n$}
We have seen that $TS^n_{\K}(F)$ is generated by $e_i^n(\ups)$ with
$\ups$ primitive. Let $R$ be generated as $\K-$algebra by
$\{r_i\}_{i\in I}$ and let  $\pi:F\to R$ be a surjective
homomorphism, then we have that the image of
\[TS_{\K}^n(\pi):TS^n_{\K}(F)\to TS^n_{\K}(R)\]
is the subalgebra of $TS_{\K}^n(R)$ generated by the
$TS^n_{\K}(\pi)(e_i^n(\ups))$ but this latter are the orbits sums of
$\pi(\ups)^{\otimes i}\otimes 1^{\otimes(n-i)}$. Hence the image of
$TS^n_{\K}(\pi)$ is the image of the canonical morphism
$\Gamma_{\K}^n(R)\to TS_{\K}^n(R)$ corresponding to the polynomial
law $r\mapsto r^{\otimes n}$. We have proven the following
\begin{proposition}
Let $F$ be free and let $F\to R$ be onto. Write
$\tau:\Gamma_{\K}(R)\to TS^n_{\K}(R)$ for the canonical homomorphism
such that $\tau(\pd r n)=r^{\otimes n}$ for all $r\in R$. Then
$\tau(\Gamma_{\K}^n(R))=TS_{\K}^n(\pi)(TS^n_{\K}(F))$.
\end{proposition}
\section{Relations: the first syzygy}
We have generators we look for relations between them: the first
syzygy.
\subsection{The free object}
When you are looking for syzygies the first step is usually to find
some free objects mapping onto your algebra. The basic idea here is
to use the basis $\mathcal{B}_n$ we found for $TS^n_{\K}(F)$
allowing $n\rightarrow \infty$.
\begin{definition}
Let $F=\K\{x_1,\dots,x_m\}$ as usual.
\begin{itemize} \item Set
$\mathfrak{F}$ for the free $\K-$module with basis
$\mathcal{B}:=\{f_{\al}\,:\, \al\in\NN^{(\up^+)}\}$.
\item We set $1_{\fr}:=f_{\mathbf{0}}$ where
$\mathbf{0}\in\NN^{\up^+}$ is such that $\mathbf{0}(\ups)=0$ for all
$\ups\in\up^+$.
\item For $i\in\NN-\{0\}$ and $\ups\in\up^+$ we denote by $f_i(\ups)$ the
element of $\mathcal{B}$ corresponding to $\al:\up^+\rightarrow \NN$
such that $\al(\ups)=i$ and $\al(\mu)=0$ otherwise.
\item We set $\fr^n$ the sub$-\K-$module generated by those $f_{\al}$
with $\sum_{\ups\in\up^+}\al(\ups)>n$.
\end{itemize}
\end{definition}
For all $n\geq 1$ there a split exact sequence of $\K-$modules
\begin{equation}\label{exmpa}
\xymatrix@C=0.5cm{
  0 \ar[r] & \fr^n \ar[rr] && \fr \ar[rr]^{\sigma_n} && TS^n_{\K}F \ar[r] & 0 }
\end{equation}
where
\[\sigma_n:\fr\rightarrow TS^n_{\K}F\] is given by
\begin{equation}
\sigma_n:\begin{cases} f_{\al}\mapsto e_{\al}^n, & \hbox{if }
\sum_{\ups\in\up^+}\al(\ups)\leq n\\
f_{\al}\mapsto 0, & \hbox{otherwise}
\end{cases}
\end{equation}
Using this splitting one can lift the product of $TS^n_{\K}F$ to
$\fr$ making it into an associative $\K-$algebra with identity
$1_{\fr}$. Observe indeed that the Product Formula stabilizes for
$n$ big enough because the number of solutions of $(\ref{sys})$ is
finite also if one drops out the constrain $\sum\gm_{ij}\leq n$.
Thus one can express $e^n_{\al}e^n_{\bt}=\sum_{\gm}e^n_{\gm}$ with
respect to $\mathcal{B}_n$ using the Product Formula with $n\gg
\max(\sum_i\al_i,\sum\bt_j)$ and then define
$f_{\al}f_{\bt}:=\sigma_n^{-1}(\sum_{\gm}e^n_{\gm})=\sum_{\gm}f_{\gm}$
\begin{proposition}
The $\K-$algebra $\fr$ is generated by $f_i(\ups)$ where $i\geq 1$
and $\ups\in\up^+$ is primitive.
\end{proposition}
\begin{proof}
By the very proof of Th.\ref{prigen} and the definition of $\fr$ it
is possible to express any $f_{\al}\in\mathcal{B}$ as an element of
the subalgebra generated by the $f_i(\ups)$ of the statement.\qed
\end{proof}
\begin{corollary}\label{gefre} The $\K-$algebra $\fr^{ab}$ is generated by
$\mathfrak{ab}(f_i(\ups))$ where $i\geq 1$ and $\ups$ that varies in
a complete set of representatives of $\Psi$.
\end{corollary}
\begin{proof}
By the above prop. and the section on generators.\qed
\end{proof}

\subsection{Degree}
Few words on degree: we have a multidegree $\partial$ on
$TS^n_{\K}(F)$ and again we can lift it up to $\fr$ using $\sigma_n$
namely by
\begin{equation}\partial(f_{\al}):=\sum_{\ups\in\up^+}\al(\ups)\partial({\ups})
\end{equation}
\begin{definition}
For $\al\in\NN^m$ write $\fr_{\al}$ (resp. $TS^n_{\K}(F)_{\al}$) for
the $\K-$submodule of $\fr$ (resp. of $TS^n_{\K}(F)$) of the
elements of multidegree $\al$. Analogously we define the total
degree of an homogeneous element of $f\in\fr$ to be the sum
$\sum_i\partial_i(f)$.
\end{definition}
By construction $$\fr=\bigoplus_{\al}\fr_{\al}$$ and this make $\fr$
a $\NN^m-$graded $\K-$algebra as well as
$$TS^n_{\K}(F)=\bigoplus_{\al}TS^n_{\K}(F)_{\al}.$$
The same holds for the abelianizations $\fr^{ab}$ and
$TS^n_{\K}(F)^{ab}$ as can be easily seen.
\begin{proposition}
One has $\fr_{\al}\cong TS^n_{\K}(F)_{\al}$ for all $\al\in\NN^m$
such that $\sum_i\al_i\leq n$ that is
$\ker\sigma_n\cap\bigoplus_{\al,\,\sum_i\al_i\leq
n}\fr_{\al}=\{0\}$. In other words $\fr_{\al}\cong
TS^n_{\K}(F)_{\al}$, $\forall n>\sum_i\al_i$.
\end{proposition}
\begin{proof}
By the definition of degree it is clear that
$\fr_{\al}\cap\fr^n=\{0\}$ for all $n>\sum_i\al_i$ the total degree
of the elements of $\fr_{\al}$ and the result follows.\qed
\end{proof}
\begin{corollary}
Let $\fr_k$ (resp. $TS^n_{\K}(F)_k$) be the linear span of the
elements of total degree $k\in\NN$ of $\fr$ (resp. of
$TS^n_{\K}(F)$), then $\fr_k\cap\fr^n=\{0\}$ for all $n>k$ and
$\fr_k\cong TS^n_{\K}(F)_k$ for all $n>k$.
\end{corollary}
\begin{proof}
Observe that $\fr_k=\bigoplus_{\al,\,\sum_i\al_i=k}\fr_{\al}$ and so
on. Then apply the above corollary.\qed
\end{proof}
\begin{proposition}\label{stru}
Let $\fr$ be endowed with the above defined product then the
sequence (\ref{exmpa}) gives an isomorphism of $\NN^m-$graded
$\K-$algebra $TS^n_{\K}(F)\cong \fr/\fr^n$. The splitting in
(\ref{exmpa}) is not an algebra homomorphism.
\end{proposition}
\begin{proof}
By construction $\sigma_n$ is a surjective $\NN^m-$graded
$\K-$algebras homomorphism whose kernel is $\fr^n$. Looking at
Product Formula one can check that $f_{\al}f_{\bt}$ does not belong
in general to the submodule generated by those $f_{\gm}$ such that
$\gm\leq\max(\sum_i\al_i,\sum\bt_j)$ hence the sequence is not split
in the category of $\K-$algebras.\qed
\end{proof}

\section{Main result}
We would like to show now that $\fr^{ab}$ is a freely generated by
the $\mathfrak{ab}(f_i(\ups))$ where $i\geq 1$ and $\ups$ that
varies in a complete set of representatives of $\Psi$. In order to
prove this result we need some machinery coming from representations
theory.
\begin{definition}
We write $f_{i,\ups}$ for $\mathfrak{ab}(f_i(\ups))$.
\end{definition}
\subsection{Generic matrices} This subsection
is borrowed from Procesi, see {\cite{dp}} for a recent paper and
{\cite{p1}} for the original sources.

Let $A_n:=\K[\xi_{ij,h}]$ be a polynomial ring where $i,j=1,\dots,n$
and $h=1,\dots,m$, note that $A_n$ is isomorphic to the symmetric
$\K$-algebra of the dual of $Mat(n,\K)^m$.

Let $F$ be again the free associative $\K$-algebra on $m$ generators
then
\[\nn_{\K}(F,Mat(n,S))\cong Mat(n,S)^m \cong
\C_{\K-alg}(A_n,S)\] for any commutative $\K$-algebra $S$. More
precisely set $B_n:=Mat(n,A_n)$ and let $\xi_h\in B_n$ be given by
$(\xi_h)_{ij}=\xi_{ij,h}$\,, $\forall\, i,j,h$, these are called the
($n\times n$) {\em generic matrices} (over $\K$) and were introduced
by Procesi (see {\cite{p1}}). Let $\pi_n:F\to B_n$ be the
$\K-$algebra given by $x_h\mapsto \xi_h$. We have then that given
any $\rho\in Hom_{\K-alg}(F,Mat(n,S)))$, with $S$ a commutative
$\K$-algebra, there is a unique $\bar{\rho}\in Hom_{\K-alg}(A_n,S)$
given by $\xi_{ij,h}\mapsto (\rho(\xi_h))_{ij}$ and it is such that
the following diagram commutes
\begin{equation}
\xymatrix{
  F \ar[dr]_{\rho} \ar[r]^{\pi}
                & B_n \ar[d]^{(\bar{\rho})_n}  \\
                & Mat(n,S)             }
\end{equation}
where $(-)_n$ denotes the induced map on $n\times n$ matrices. The
homomorphism $\pi_n$ is called the universal $n$ dimensional
representation (for the free algebra). We denote by $\mathcal{G}_n$
the subring of $B_n$ generated by the generic matrices i.e. the
image of $\pi_n$.
\subsection{Determinant} The composition $det\cdot\pi$ gives a multiplicative polynomial law $F\to A$
homogeneous of degree $n$ hence a unique homomorphism
$\delta_n:TS^n_{\K}(F)^{ab}\to A$ such that
$\delta_n(\mathfrak{ab}(e_n^n(f)))=det(\pi_n(f))$.
\begin{definition} Let $C_n\subset A_n$ be the subalgebra generated by
the coefficients of the characteristic polynomials of the elements
of $\mathcal{G}_n$. We write
\[det(t-f)=t^n+\sum_{i=1}^n(-1)^i\varsigma_i^n(f)t^{n-i}\]
where $f\in\mathcal{G}_n$.
\end{definition}
Note that
$\delta_n(\mathfrak{ab}(e_i^n(f)))=\varsigma_i^n(\pi_n(f))$ for all
$f\in F$.

By Theorem \ref{genab} we have that
$\delta_n(TS^n_{\K}(F)^{ab})=C_n$. If one gives $\xi_{ij,h}$
multidegree $\partial(\xi_{ij,h})=\partial(x_h)$ then
$A_n=\bigoplus_{\al\in\NN^m} A_{n,\al}$ is a $\NN^n-$graded ring
with homogeneous components $A_{n,\al}$ and one can easily check
that $C_n=\bigoplus_{\al\in\NN^m} C_{n,\al}$ where
$C_{n,\al}=A_{n,\al}\cap C_n$ are the homogeneous component.

One checks that $\delta_n(TS^n_{\K}(F)_{\al}^{ab})=C_{n,\al}$.

For all $n$ there is an homomorphism of $\NN^m-$graded rings
$\varphi_n:A_n\twoheadrightarrow A_{n-1}$ given by
$x_{nj,h},x_{in,h}\mapsto 0$ and $x_{ij,h}\mapsto x_{ij_h}$ for
$i,j< n$.
\begin{definition}
We write $C_{\al}:=\lim_n C_{n,\al}$ and
$C=\bigoplus_{\al\in\NN^m}C_{\al}$ and
$\varsigma_i(f)=\lim_n\varsigma_i^n(f)$.
\end{definition}
\begin{remark}
The algebra $C$ is the graded inverse limit of the $C_n$.
\end{remark}
\begin{proposition}\label{donk}
The ring $C$ is a polynomial ring freely generated by the
$\varsigma_{i,\ups}$ with $i\geq 1$ and $\ups$ that varies in a
complete set of representatives of $\Psi$.
\end{proposition}
\begin{proof}
By {\cite{d2}},\S3,(10)and remark loc.cit. and
Complements{\cite{d1}}.\qed
\end{proof}
\subsection{Freeness} We are now in the following situation: there are surjective
homomorphisms
\[\fr^{ab}\xrightarrow{\tilde{\sigma_n}}
TS^n_{\K}(F)^{ab} \xrightarrow{\delta_n} C_n\] of $\NN^m-$graded
algebras that commutes with $\varphi_n$ hence there is an
homomorphism $\delta:\fr^{ab}\to C$ such that
$f_{i,\ups}\mapsto\varsigma_i(\ups)$ and this is an isomorphism onto
$C$ with inverse given by $\varsigma_i(\ups)\mapsto f_{i,\ups}$ by
Prop.\ref{donk} and Cor. \ref{gefre}. We have proven the following
\begin{proposition}
The algebra $\fr^{ab}$ is a free polynomial ring freely generated by
$f_{i,\ups}$ with $i\geq 1$ and $\ups$ that varies in a complete set
of representatives of $\Psi$.
\end{proposition}
\begin{remark}
The above Proposition can be also proved observing that
$\fr\cong\hat{\Gamma}(F^+)$ where the latter has been introduced by
Ziplies {\cite{zipgen}} and applying then Th.4.4{\cite{zipab}}.
\end{remark}
\subsection{A presentation}
The assignment $R\mapsto R^{ab}$ gives a functor from $\K-$algebras
to commutative $\K-$algebras. It preserves surjections and for
$\rho:R\rightarrow S$ we denote by $\tilde{\rho}:R^{ab}\rightarrow
S^{ab}$ its image by $ab$. Recall that we denote by
$\mathfrak{ab}:R\rightarrow R^{ab}$ the canonical projection. We
finally state the main result of this section
\begin{theorem}\label{prese}
Let $\K$ be a commutative ring and let $F$ be a free $\K-$algebra
then
\begin{equation}
\xymatrix@C=0.5cm{
  0 \ar[r] & (\fr^n)^{ab} \ar[rr] && (\fr)^{ab} \ar[rr]^{\tilde{\sigma}_n} && TS^n_{\K}F^{ab} \ar[r] & 0 }
\end{equation}
is a presentation by generators and relations.
\end{theorem}
\begin{corollary}
Let $R\cong F/I$ be a $\K-$algebra generated by $\{r_j\}$ then we
have a presentation by generators and relations
\begin{equation}
\xymatrix@C=0.5cm{
  0 \ar[r] & (\fr^n)^{ab}+J \ar[rr] && (\fr)^{ab} \ar[rr] && \Gamma_{\K}^n(R)^{ab} \ar[r] & 0 }
\end{equation}
where $J$ is linearly generated by the lifting of
$e_{\al}^n(f_1,\dots,f_k)$ where at least one $f_h\in I$.
\end{corollary}
\begin{proof}
By {\cite{bo}} es.2)(f) pag.92.\qed
\end{proof}
\begin{remark}
the above corollary gives generators and relations for
$\Gamma_{\K}^n(R)$ for any commutative ring $R$.
\end{remark}
\begin{corollary}\label{fields}
Let $\K$ be an infinite field then $\fr^n$ and $(\fr^n)^{ab}$ are
generated as ideals by the lifting of the $e_{n+k}^h(f)$ with $f\in
F$ and their images by $\mathfrak{ab}$ respectively.
\end{corollary}
\begin{proof}
Over an infinite field one has that $TS_{\K}^n(M)$ is linearly
generated by $m^{\otimes n}$ as $m\in M$ for any $M\in
Mod_{\K}$.\qed
\end{proof}
\section{Invariants}
The general linear group $G:=GL(n,\K)$ made of the invertible
matrices of $Mat(n,\K)$ acts on $Mat(n,\K)^m$ by simultaneous
conjugation, i.e. via base change on $\K^n$. Recall that the
categorical quotient $\M_n^m//G$ is defined as
\begin{equation}\M_n^m//G:=Spec(A^G)\end{equation}
where as usual $A_n^G$ denotes the ring of the invariants,
$\M_n^m=Spec(A_n^G)$ and the action of $G$ on $A_n$ is induced by
the one on $Mat(n,\K)^m$. This scheme has be widely studied and we
refer the reader to {\cite{ar,d1,p1}} for masterpieces on this
subject. It is the coarse moduli space parameterizing the
$n-$dimensional linear representations of $F$ up to base change.

Being the determinant invariant under base change we have that the
ring $C_n$ is made of invariants i.e $C_n\subset A_n^G$. When $\K$
is a characteristic zero field was showed firstly by C. Procesi
{\cite{p1}} and separately Sibirsky that $C_n=A_n^G$. This has been
generalized to the positive characteristic case by S. Donkin and
then by A. Zubkov proving a Procesi's conjecture, see
{\cite{d1,p1,zucom}}. It also prove that the result remain true over
the integers showing indeed that $C_n$ is a $\mathbb{Z}$ form of the
ring of invariants, i.e. $\K\otimes_{\mathbb{Z}}C_n\cong
\K[Mat(n,\K)^m]^{GL(n,\K)}$. Zubkov in \cite{zu}, gave a
generalization of the well celebrated Procesi-Razmyslov theorem
{\cite{pa,ra}}, namely he proved that the kernel of the surjection
$C\to C_n\cong A_n^G$ is generated by $\psi_{n+k}(f)$ with $k\geq 1$
when $\K$ is an infinite field.

Combining these results on the invariants of matrices with Theorem
\ref{prese} and Corollary\ref{fields} we finally have

\begin{theorem}\label{end}
Let $\K$ be an infinite field or the ring $\Z$ of integers. The
homomorphism $TS^n_{\K}(F)^{ab}\to \K[Mat(n,\K)]^{GL(n,\K)}$ induced
by the composition $det\cdot \pi_n$ of the determinant with the
universal representation is an isomorphism
\[TS^n_{\K}(F)^{ab}\xrightarrow{\cong} \K[Mat(n,\K)]^{GL(n,\K)}\].
\end{theorem}
\begin{proof}
When $\K$ is an infinite field it follows by the results on the
invariants of matrices recalled above, Theorem \ref{prese} and
Corollary\ref{fields}. It remains to prove only the case of $\K=\Z$.
For any commutative ring $\K$ one has
$TS^n_{\K}(F)\cong\K\otimes_{\Z}TS^n_{\Z}(F)$ then the homomorphism
$id_{\K}\otimes \mathfrak{ab}_{\Z}:\K\otimes_{\Z}
TS^n_{\Z}(F)\to\K\otimes_{\Z} TS^n_{\Z}(F)^{ab}$ factors through
$\mathfrak{ab}_{\K}:TS^n_{\K}(F)\to TS^n_{\K}(F)^{ab}$ by
universality of $ab$ and the following diagram commutes
\begin{equation}
\xymatrix{
  TS^n_{\K}(F) \ar[d]_{\mathfrak{ab}_{\K}} \ar[r]^(.4){\cong}
                & **[rr]\K\otimes_{\Z} TS^n_{\Z}(F) \ar[d]^{id_{\K}\otimes \mathfrak{ab}_{\Z}}  \\
  TS^n_{\K}(F)^{ab} \ar[r]_(.4){\varepsilon_{\K}}
                & **[rr]\K\otimes_{\Z} TS^n_{\Z}(F)^{ab}             }
\end{equation}
Let $\K$ be an algebraically closed field.

We have a map $TS^n_{\Z}(F)^{ab}\to\Z[Mat(n,\Z)^m]^{GL(n,\Z)}$
induced by the composition of the universal $n-$dimensional
representation with the determinant, this gives another one
$\K\otimes_{\mathbb{Z}}TS^n_{\Z}(F)^{ab}\to\K\otimes_{\mathbb{Z}}\Z[Mat(n,\Z)^m]^{GL(n,\Z)}$
that can be lifted to
$TS^n_{\K}(F)\cong\K[Mat(n,\K)^m]^{GL(n,\K)}\cong\K\otimes_{\mathbb{Z}}\Z[Mat(n,\Z)^m]^{GL(n,\Z)}$
by the previous discussion and {\cite{d1}}Compl. We have then
$\K\otimes_{\mathbb{Z}}TS^n_{\Z}(F)^{ab}\cong\K\otimes_{\mathbb{Z}}\Z[Mat(n,\Z)^m]^{GL(n,\Z)}$
for any algebraically closed field as $\NN^m-$graded rings with
finitely generated homogenous summand and the result follows.\qed
\end{proof}

\section{acknowledgements} I would like to thank C. Procesi, M. Brion
and C. De Concini for hints and useful discussions. A special thank
to David Rydh for valuable comments on the first version of this
paper.

\end{document}